\documentclass[11pt]  {amsart} 
\usepackage{amsmath,latexsym,amssymb,amsmath, 
amscd,amsthm,amsxtra}

\newtheorem{theorem}{Theorem}[section]

\newtheorem{remark}{\textbf{Remark}}[section]
\def\o{\omega}
\def\O{\Omega}
\def\p{\partial}

\def\a{\alpha}

\def\De{\Delta}
\def\n{\nabla}

\def\<{\langle}
\def\>{\rangle}
\def\div{{\rm div}}
\def\De{\Delta}
\def\n{\nabla}
\def\ode{\overline{\Delta}}
\def\on{\overline{\nabla}}

\def\RR{\mathbb{R}}

\def\SS{\mathbb{S}}
\def\HH{\mathbb{H}}

\def\o{\omega}
\def\O{\Omega}
\def\p{\partial}

\def\a{\alpha}

\begin{document}

\title{A Minkowski type inequality in space forms}

\author{Chao Xia}\address{School of Mathematical Sciences, Xiamen University, 361005, Xiamen, China and
 Department of Mathematics and Statistics, McGill University, Montreal, H3A 0B9, Canada}
 \email{chaoxia@xmu.edu.cn}
\thanks{Research of CX is  supported in part by the Fundamental Research Funds for the Central Universities (Grant No. 20720150012), NSFC (Grant No. 11501480)  and CRC Postdoc Fellowship. }

\begin{abstract}
In this note we apply the general Reilly formula established in \cite{QX} to the solution of a Neumann boundary value problem to prove an optimal Minkowski type inequality in space forms.
\end{abstract}

\date{}
\maketitle

\medskip


\section{Introduction}
  Let $(\Omega^{n}, g)$ be an $n$-dimensional compact  Riemannian manifold with smooth boundary $\p\O= M$. Let $H$ be the (normalized) mean curvature  and $h$ be the second fundamental form of $M\subset \O$ respectively. In the paper \cite{QX}, we (joint with Qiu) have proved the following generalization of Reilly's formula. We use the same notations as in \cite{QX}. 
  
\vspace{3mm}
  
\noindent{\bf Theorem A.} (Qiu-Xia \cite{QX}) {\it Let $V: \O\to \RR$ be a given a.e. twice differentiable function. Given a smooth function $f$  on $\O$,  we denote $z=f|_M$ and $u=\on_\nu f$. Let $K\in\RR$. Then we have the following identity:
\begin{eqnarray}\label{qx}
&&\int_\O V\left((\ode f+Knf)^2-|\on^2 f+Kfg|^2\right)d\O\nonumber\\&=&\int_M V\left(2u\De z+(n-1)Hu^2+h(\n z, \n z)+(2n-2)Kuz\right)dA\nonumber\\&&+\int_M \on_\nu V\left(|\n z|^2-(n-1)Kz^2\right) dA\nonumber\\&&+\int_\O \left(\on^2 V-\ode Vg-(2n-2)KVg+V Ric\right)(\on f, \on f)d\O\nonumber\\&&+(n-1)\int_\O (K\ode V+nK^2 V)f^2 d\O.
\end{eqnarray}
}

When $V\equiv 1$ and $K=0$, \eqref{qx} reduces to the classical Reilly's formula \cite{Reilly0, Reilly1}. Reilly's original formula  has numerous applications, see for example \cite{Reilly1, CW, PRS, Ros, Reilly2}. In \cite{QX}, we successfully apply the general Reilly formula \eqref{qx} to prove a new Heintze-Karcher type inequality for compact manifolds with mean convex boundary and sectional curvature bounded below.  In this paper, we continue to explore other applications of \eqref{qx}.

Reilly \cite{Reilly2} used his formula to prove the following Minkowski inequality for compact Riemannian manifolds with nonnegative Ricci curvature and convex boundary.

\vspace{3mm}
\noindent{\bf Theorem B.} (Reilly \cite{Reilly2}) {\it Let $(\Omega^{n}, g)$ be a compact $n$-dimensional Riemannian manifold with smooth convex boundary $M$ and non-negative Ricci curvature. 
Then
\begin{eqnarray}\label{rm}
{\rm Area}(M)^2\geq n{\rm Vol} (\Omega) \int_M H dA.
\end{eqnarray}
 The equality in \eqref{rm} holds if and only if $\O$ is isometric to an Euclidean ball.
}

\vspace{3mm}
When $\O\subset\mathbb{R}^n$, inequality \eqref{rm} is exactly a special case of Minkowski's inequality for mixed volumes in the theory of convex bodies, see \cite{Sch}, Theorem 7.2.1. A diffenrent proof of Theorem B was given by Wang-Zhang \cite{WZ}, based on the Alexandrov-Bakelman-Pucci estimate.

Reilly's proof  is based on the solvability of the following Neumann problem 
\begin{equation}\label{Neumann0}
\left\{
\begin{array}{rccl}
\ode f&=&1&
\hbox{ in } \Omega,\\
u&=&c &\hbox{ on } \partial \Omega,\\
\end{array}
\right.
\end{equation}
for $c=\frac{{\rm Vol}(\O)}{{\rm Area}(\p\O)}$.
He applied his formula \eqref{qx} ( for $K=0$ and $V\equiv 1$) to the solution  of \eqref{Neumann0} to derive
$$\frac{n-1}{n}{\rm Vol}(\O)\geq c^2\int_M HdA,$$ which is \eqref{rm}.

The topic of geometric inequalities for curvature integrals in non-Euclidean space forms attracts many attentions in recent years, see for example \cite{WX} and refenreces therein.
Curvature integral with ``weight" seems quite natural in the general relativity, especially in the hyperbolic space. Quite recently, Brendle-Hung-Wang \cite{BHW} established a Minkowski type inequality between ``weighted" mean curvature integral and ``weighted" volume for hypersurfaces in anti-de Sitter-Schwarzschild manifolds by using a ``weighted" Heintze-Karcher inequality by Brendle \cite{Brendle}.  See also \cite{dLG, GWW, GWWX, QX} for related works.

In this shote note, based on Theorem A, we  prove an analog of Minkowski's inequality for ``weighted mixed volumes" in non-Euclidean space forms. We use $\HH^n$ to denote the hyperbolic space with curvature $-1$ and $\SS_+^n$ to denote the open hemi-sphere with curvature $1$.
\begin{theorem}\label{rm22} Let $\Omega^{n}\subset \mathbb{H}^n$ ( $\SS_+^n$ resp.) be a compact $n$-dimensional domain with smooth boundary $M$. 
Let $V(x)=\cosh r$ ($\cos r$ resp.), where $r(x)=dist(x,p)$ for some fixed point $p\in\HH^n$ ($p\in{\mathbb{S}_+^n}$ resp.). We further assume the second fundamental form of $M$ satisfies \begin{eqnarray}\label{ass}
h_{ij}\geq \on_\nu \log V g_{ij}.
\end{eqnarray}
Then we have
\begin{eqnarray}\label{rm2}
\left(\int_M VdA\right)^2\geq n\int_\Omega Vd\O \int_M HV dA.
\end{eqnarray}
The equality in \eqref{rm2} holds if and only if $\O$ is a geodesic ball $B_R(q)$ for some point $q\in \HH^n$ ($q\in{\mathbb{S}_+^n}$ resp.). In particular,   \eqref{rm2} holds true when $M$ is horo-spherical convex in the case $\O\subset\HH^n$ or $M$ is convex and $p\in \O$ in the case $\O\subset\SS_+^n$.
\end{theorem}

The horo-spherical convexity of $M\subset \HH^n$ means that all the principal curvatures are bigger than or equal to $1$. Condition \eqref{ass} seems like some kind of convexity for $M$. Particularly, when $\O\subset \mathbb{R}^n$ and $V\equiv 1$, this is the usual convexity. Moreover,  horo-convexity in $\HH^n$ and convexity in $\SS_+^n$ imply  condition \eqref{ass}. This follows because $\on_\nu V<V$ in the case $\O\subset\HH^n$ and $\on_\nu V\leq 0$ in the case $p\in\O\subset \SS_+^n$. 
We remark that, the equality in \eqref{rm2} holds for not only  geodesic balls centered at $p$ but all geodesic balls.

In the Euclidean space, Theorem B is equivalent to say that \begin{eqnarray}\label{concave}
\frac{d^2}{dt^2}{\rm Vol}(\O_t)^{\frac1n}\leq 0,
\end{eqnarray}
where $\O_t=\O+tB=\{x\in \mathbb{R}^n| dist(x,\O)\leq t\}.$
 Similarly, Theorem 1 can be interpreted as the following equivalent statement.

\begin{theorem}\label{rm33}
Let $\Omega^{n}\subset \mathbb{H}^n$ ( $\SS_+^n$ resp.) and $V$ be as in Theorem \ref{rm22}. Let $K=-1$ ($K=1$ resp.). Denote  $\Omega_t:=\{x\in \mathbb{H}^n (\SS_+^n \hbox{ resp.)} | dist(x, \O)\leq t\}.$ For the case $\SS_+^n$ we assume $t\in [0,T)$ for which $\Omega_t\subset \SS_+^n$. Then 
\begin{eqnarray*}
\frac{d^2}{dt^2}\left(\int_{\O_t} V d\O\right)^{\frac1n}+K\left(\int_{\O_t} V d\O\right)^{\frac1n}\leq 0.
\end{eqnarray*}
\end{theorem}

\

The idea to prove \eqref{rm2} is parallel to Reilly's. We will ultilize the solution to a Neumann boundary value problem \eqref{Neumann} and the general Reilly formula. However, the computation is much more complicated due to the complication of the boundary terms in the general Reilly formula.

\


\section{Proof of Theorem \ref{rm22} and \ref{rm33}}
Let $V=\cosh r, K=-1$ or $V=\cos r, K=1$ in \eqref{qx} for the case $\HH^n$ or $\SS_+^n$ respectively, where $r(x)= dist(x,p)$. The function $f$ is the solution to the following Neumann boundary value problem:
\begin{equation}\label{Neumann}
\left\{
\begin{array}{rccl}
\ode f+Knf&=&1&
\hbox{ in } \Omega,\\
Vf_\nu-V_\nu f&=&cV  &\hbox{ on } \partial \Omega,\\
\end{array}
\right.
\end{equation}
where $V_\nu:=\on_\nu V$ and $c=\frac{\int_{\O} V}{\int_M V}.$
We claim that there exists a unique solution $f\in C^\infty(\O)$ to \eqref{Neumann}, up to an additive $\a V$ for constants $\a\in \mathbb{R}$.
In fact, it follows from the Fredholm alternative that there exists a unique solution $w\in C^\infty(\O)$ (up to an additive constant) to the following Neumann boundary value problem
\begin{equation}\label{Neumann1}
\left\{
\begin{array}{rccl}
\div( V^2\on w)&=&V&
\hbox{ in } \Omega,\\
V^2w_\nu&=&cV  &\hbox{ on } \partial \Omega.\\
\end{array}
\right.
\end{equation}if and only if $c=\frac{\int_{\O} V}{\int_M V}$. Using \eqref{vv} below, one checks readily that $f=wV$ solves \eqref{Neumann}.

 For simplicity, we omit the volume form $d\O$ and the area form $dA$ in the intergrations.

It is well-known that $V$ satisfies
\begin{eqnarray}\label{vv}
\on^2 V=-KVg, 
\end{eqnarray}
which will be used frequently in the following.

We will use the solution $f$ of $\eqref{Neumann}$ in the general Reilly formula \eqref{qx}.
For our choice of $K$ and $V$, we see from \eqref{vv} that the integrand in last two lines of \eqref{qx} vanishes.
By using H\"older's inequality and the equation in \eqref{Neumann}, we have from \eqref{qx} that
\begin{eqnarray}\label{R1}
\frac{n-1}{n}\int_\O V&\geq &\int_\O V\left((\ode f+Knf)^2-|\on^2 f+Kfg|^2\right)\nonumber\\&=&\int_{M} V \left(2u\De z+(n-1)Hu^2+h(\n z, \n z)+(2n-2)Kuz\right)\nonumber\\&&+\int_M V_\nu\left(|\n z|^2-(n-1)Kz^2\right).
\end{eqnarray}

Let us investigate the RHS of \eqref{R1}. 
By using the Gauss-Weigarten formula and \eqref{vv}, we see \begin{eqnarray}\label{ss}
&&\n_i V_\nu=\on_i\on_\nu V+h_{ij}V_j=h_{ij}V_j,\end{eqnarray}
\begin{eqnarray}\label{rr}
&&\De V=\ode V- \on_\nu\on_\nu V-(n-1)HV_\nu=-(n-1)KV-(n-1)HV_\nu.
\end{eqnarray}

Using the Neumann boundary condition in \eqref{Neumann}, integration by parts, \eqref{ss} and \eqref{rr}, we have
\begin{eqnarray}\label{R2}
&&\int_{M} 2Vu\De z=\int_M 2(V_\nu z+cV)\De z\nonumber\\&=&\int_M -2V_\nu |\n z|^2-2z\n V_\nu \n z+2cz\De V\nonumber\\&=&\int_M -2V_\nu |\n z|^2-2zh_{ij}V_i z_j+2(n-1)cz(-KV-HV_\nu),
\end{eqnarray}
\begin{eqnarray}\label{R3}
&&\int_{M} (n-1)Hu^2V=\int_M (n-1)HV\left(c+\frac{V_\nu}{V}z\right)^2\nonumber\\&=&\int_M (n-1)c^2HV+2c(n-1)HV_\nu z+(n-1)H\frac{V_\nu^2}{V}z^2,
\end{eqnarray}
\begin{eqnarray}\label{R4}
&&\int_{M} (2n-2)KuzV=\int_M 2(n-1)Kz\left(cV+V_\nu z\right).
\end{eqnarray}

Inserting \eqref{R2}-\eqref{R4} into \eqref{R1}, we have
\begin{eqnarray}\label{R5}
\frac{n-1}{n}\int_\O V&\geq &\int_{M} -V_\nu |\n z|^2-2zh_{ij}V_i z_j+ (n-1)c^2HV\nonumber\\&&+(n-1)H\frac{V_\nu^2}{V}z^2+h(\n z, \n z)V+(n-1)KV_\nu z^2.
\end{eqnarray}

 Multiplying $-\frac{V_\nu}{V}z^2$ to both side of \eqref{rr}, integrating by parts and using \eqref{ss}, we have
 \begin{eqnarray}\label{R6}
&&\int_M (n-1)H\frac{V_\nu^2}{V}z^2 +(n-1)KV_\nu z^2\nonumber\\&=&\int_M -\frac{V_\nu}{V}z^2\big[-(n-1)KV-(n-1)HV_\nu\big]\nonumber\\&=& \int_M  -\frac{V_\nu}{V}z^2 \De V\nonumber\\&=&\int_M  \frac{\n V_\nu \n V}{V}z^2+ \frac{2z\n z\n V V_\nu}{V}-\frac{V_\nu z^2}{V^2}|\n V|^2\nonumber\\&=&\int_M  \frac{h_{ij}V_iV_j}{V}z^2+ \frac{2z\n z\n V V_\nu}{V}-\frac{V_\nu z^2}{V^2}|\n V|^2.
\end{eqnarray}

Inserting \eqref{R6} into \eqref{R5}, we obtain
\begin{eqnarray}\label{R7}
\frac{n-1}{n}\int_\O V&\geq &\int_{M} -V_\nu |\n z|^2-2zh_{ij}V_i z_j+ (n-1)c^2HV\nonumber\\&&+h(\n z, \n z)V+\frac{h_{ij}V_iV_j}{V}z^2+ \frac{2z\n z\n V V_\nu}{V}-\frac{V_\nu z^2}{V^2}|\n V|^2\nonumber\\&=&\int_{M} (n-1)c^2HV \nonumber\\&&+Vh_{ij}\left(z_i-\frac{V_iz}{V}\right)\left(z_j-\frac{V_jz}{V}\right)-V_\nu\left|\n z-\frac{\n Vz}{V}\right|^2.
\end{eqnarray}
By the assumption \eqref{ass}, the last line in \eqref{R7} is nonnegative.
Therefore, we derive from \eqref{R7} that
\begin{eqnarray}\label{R8}
\frac{n-1}{n}\int_\O V\geq \int_{M} (n-1)c^2HV=\frac{\left(\int_\O V\right)^2}{\left(\int_M V\right)^2} \int_{M} (n-1)HV.
\end{eqnarray}
It follows that
\begin{eqnarray}\label{R9}
\left(\int_M V dA \right)^2\geq n\int_\O V d\O\int_{M} HV dA.
\end{eqnarray}

Let us explore the equality case in \eqref{R9}. We consider the case $\O\subset\HH^n$. First, for a geodesic ball $B_R(p)\subset \HH^n$, centered at $p$, $V=\cosh R$ and $H=\coth R$ are  constants on $\p B_R(p)$.
Thus $\int_{\p B_R(p)} V dA=\o_{n-1}\cosh R\sinh^{n-1} R$ and $\int_{M} HV dA=\o_{n-1}\cosh^2 R\sinh^{n-2} R$. On the other hand,
$$\int_{B_R(p)} \cosh r(x) d\O(x)= \int_0^R \o_{n-1}\cosh \rho\sinh^{n-1}\rho d\rho =\frac{\o_{n-1}}{n}\sinh^n R.$$Thus equality in \eqref{R9} holds when $\O=B_R(p)$. Second, for a geodesic ball $B_R(q)\subset \HH^n$, centered at $q\in \HH^n$, not necessarily $p$, $H=\coth R$ is constant on $\p B_R(q)$ while $V$ is not. Nevertheless, we still have the quality. Indeed, by Minkowski's formula and the constancy of $H$,
$$\int_{\p B_R(q)} VdA=\int_{\p B_R(q)} HV_\nu dA= H\int_{\p B_R(q)} V_\nu dA=nH\int_{ B_R(q)} Vd\O.$$ 
Thus \begin{eqnarray*}
&&\left(\int_{\p B_R(q)} VdA\right)^2= nH\int_{ B_R(q)} Vd\O \int_{\p B_R(q)} VdA=n\int_{ B_R(q)} Vd\O \int_{\p B_R(q)} HVdA.
\end{eqnarray*}

Conversely, if the equality in \eqref{R9} holds, then by checking the equality in \eqref{R5} and \eqref{R8} we see
\begin{equation*}
\left\{
\begin{array}{rccl}
\on^2_{ij} f-fg_{ij}&=&\frac1n g_{ij}&
\hbox{ in } \Omega,\\
 \n z-\frac{z\n V}{V}&=&0  &\hbox{ on } \partial \Omega.\\
\end{array}
\right.
\end{equation*}
The boundary identity means $z=\a V$ for some constant $\a\in \mathbb{R}$.
Thus the function $\tilde{f}=f-\a V+\frac1n$ satifies
\begin{equation*}
\left\{
\begin{array}{rccl}
\on^2_{ij} \tilde{f}-\tilde{f}g_{ij}&=&0&
\hbox{ in } \Omega,\\
\tilde{f}|_{\p\O}&=&\frac1n  &\hbox{ on } \partial \Omega.\\
\end{array}
\right.
\end{equation*}
It follows from an Obata type result (see Reilly \cite{Reilly2}) that $\O$ must be some geodesic ball.

The case $\O\subset\SS^n_+$ is similar. We finish the proof of Theorem \ref{rm22}.
\qed

\

\noindent{\it Proof of Theorem \ref{rm33}:}
$\O_t$ can be viewed as the normal flow
$$\p_t X(x,t)=\nu(x,t), \quad x\in \p\O.$$ 
The variational formulas give
\begin{eqnarray}\label{var1}
\frac{d}{dt} \int_{\O_t} V d\O=\int_{\p\O} VdA_t,
\end{eqnarray}
\begin{eqnarray}\label{var2}
\frac{d}{dt}\int_{\p\O} VdA_t\nonumber&=&\int_{\p\O} V_\nu+(n-1)HVdA_t \\&=& \int_{\p\O} (n-1)HVdA_t-nK\int_{\O_t} Vd\O.
\end{eqnarray}
Using \eqref{var1}, \eqref{var2} and Theorem \ref{rm2}, we deduce
\begin{eqnarray}
\frac{d^2}{dt^2} \left(\int_{\O_t} V d\O\right)^{\frac1n}&=& \frac1n\left(\int_{\O_t} V d\O\right)^{\frac1n-1} \left(\int_{\p\O} (n-1)HVdA_t -nK\int_{\O_t} V d\O\right)\nonumber\\&&+\frac1n(\frac1n-1)\left(\int_{\O_t} V d\O\right)^{\frac1n-2} \left(\int_{\p\O_t} V dA_t\right)^2
\nonumber\\&\leq &-K\left(\int_{\O_t} V d\O\right)^{\frac1n}.\nonumber
\end{eqnarray}
We complete the proof.

\qed

\begin{remark}
It is well known that one may derive the isoperimetric inequality from \eqref{concave} in the Euclidean space. Indeed,
using the ODE comparison, we obtain
\begin{eqnarray}\label{isop}
{\rm Vol}(\O_t)^\frac1n\leq {\rm Vol}(\O)^\frac1n+\frac1n {\rm Vol}(\O)^{\frac1n-1}Area(\p\O)t.
\end{eqnarray}
Dividing both sides of \eqref{isop} by $t$ and letting $t\to\infty$, we obtain
\begin{eqnarray*}
\frac1n {\rm Vol}(\O)^{\frac1n-1}Area(\p\O)\geq\lim_{t\to\infty}\frac{1}{t}{\rm Vol}(\O_t)^\frac1n=Vol(B)^\frac1n,
\end{eqnarray*}
which is the classical isoperimetric inequality.  Similarly, in $\HH^n$, we have
\begin{eqnarray}\label{isopH}
&&\left(\int_{\O_t} V d\O\right)^{\frac1n}\leq  \left(\int_{\O} V d\O\right)^\frac1n \cosh t+ \left(\frac1n \left(\int_{\O} V d\O\right)^{\frac1n-1}\int_{\p\O} V dA\right) \sinh t.
\end{eqnarray}
However, we are not able to derive an optimal inequality between $\int_{\p\O} V dA$ and $\int_{\O} V d\O$ from \eqref{isopH} because $\lim_{t\to\infty}\frac{1}{\sinh t}\left(\int_{\O_t} V d\O\right)^{\frac1n}$ is not a dimensional constant in this case. In a forthcoming paper, we will use the flow approach to establish such kind of optimal inequality.
\end{remark}
\

{\bf Acknowledgements:} I would like to thank Professor Pengfei Guan for stimulating discussions and Prof. Guofang Wang for useful comments and for their constant supports.

\end{document}